\documentclass[12pt]{article}
\usepackage{amsfonts,amsmath,amssymb,amsthm} 

\theoremstyle{plain}
\newtheorem{thm}{Theorem}[section]
\newtheorem{lem}{Lemma}[section]
\newtheorem{cor}{Corollary}[section]
\theoremstyle{remark}
\newtheorem{rem}{Remark}[section]
\newtheorem{exmp}{Example}[section]

\def\P{{\mathbb {P}}}                       
\def\E{{\mathbb {E}}}                       \def\D{{\mathbb {D}}}
                       
\def\N{{\mathbb {N}}}                       \def\R{{\mathbb {R}}}
\def\I{{\mathbb {I}}}                       

 \def\O{{\rm{O}}}
 \def\o{{\rm{o}}}

\def\AD{{\mathcal A}}

\def\CC{{\rm\kern.24e  \vrule width.02em  height1.4ex depth-.05ex \kern-.26em C}}
\def\QQ{{\rm\kern.24em \vrule width.02em  height1.4ex depth-.05ex \kern-.26em Q}}
\def\PP{{\rm I\kern-.25em P}}         \def\RR{{\rm I\kern-.25em R}}
\def\DD{{\rm I\kern-.25em D}}         \def\EE{{\rm I\kern-.25em E}}
\def\FF{{\rm I\kern-.25em F}}         \def\NN{{\rm I\kern-.23em N}}
\def\RRp{{\rm I\kern-.25em R}_{+}}
\def\IND{{\rm 1\kern-.25em I}}

\def\etab{{\boldsymbol {\eta}}}   
\def\alphab{{\boldsymbol {\alpha}}}  
\def\kb{{\boldsymbol {k}}} 
\def\sb{{\boldsymbol {s}}} 
\def\nb{{\boldsymbol {n}}}

\begin{document}

\baselineskip 5.5truemm

\vskip0.4truecm \centerline{\bf STRONG LIMIT THEOREMS IN THE MULTI-COLOUR} 
\centerline{\bf GENERALIZED ALLOCATION SCHEME}

\vskip0.5truecm

\centerline{
  {\sc Alexey Chuprunov} and {\sc Istv\'an Fazekas}}

\bigskip
Department of Math. Stat. and Probability, Chebotarev Inst. of Mathematics and Mechanics,
Kazan State University, Universitetskaya 17, 420008 Kazan, Russia, e-mail: achuprunov@mail.ru

Faculty of Informatics, University of Debrecen,  P.O. Box 12, 4010 Debrecen, Hungary,
e-mail: fazekasi@inf.unideb.hu, tel: 36-52-512900/75211

\bigskip

\begin{abstract}
The generalized allocation scheme is studied.
Its extension for coloured balls is defined.
Some analogues of the Law of the Iterated Logarithm and the Strong Law of Large Numbers are obtained for the number of boxes containing fixed numbers of balls.
\end{abstract}

\renewcommand{\thefootnote}{}
{\footnotetext{


The publication was supported by the T\'AMOP-4.2.2.C-11/1/KONV-2012-0001 project. The project has been supported by the European Union, co-financed by the European Social Fund.
}}

\section{Introduction}
\label{sectIntr} \setcounter{equation}{0} 
The generalized allocation scheme is widely studied (see \cite{Kol99}, \cite{Kol03}, \cite{pavlov}).
The scheme contains several special cases such as the usual allocation scheme (see \cite{Weiss}, \cite{Renyi}, \cite{KSC})
and the random forests (see \cite{Kol99}, \cite{pavlov}).

In the generalized scheme of allocations of balls into boxes, the distribution of the contents of the boxes is represented as the conditional distribution
of independent random variables given that their sum is fixed. 
More precisely, the generalized allocation scheme is defined as follows. 
Let $\eta_1, \dots,\eta_N$ be non-negative integer-valued random variables. 
If there exist independent random variables $\xi_1,\dots,  \xi_N$ such that the  joint  distribution of $\eta_1, \dots, \eta_N$ admits the representation
\begin{equation}  \label{EtaXi0}
\P\{\eta_1=k_1, \dots, \eta_N=k_N\}  =  \P\{\xi_1=k_1, \dots, \xi_N=k_N\,\,\, |\,\,\, \xi_1+\dots +\xi_N=n\},
\end{equation}
where $k_1,\dots, k_N$ are arbitrary non-negative integers with $k_1+\dots + k_N=n$, we say that the distribution of $\eta_1, \dots,\eta_N$  is represented by a generalized  allocation
scheme with parameters $n$ and $N$, and independent random variables $\xi_1,\dots, \xi_N$. 
In view of independence of the random variables $\xi_1,\dots, \xi_N$, the study of several
questions of the generalized allocation scheme can be reduced to problems of sums of independent random variables. 
Let $\mu_{snN}$ be the number of the random variables $\eta_1, \dots, \eta_N$  being equal to $s$, $s=0,1,\dots, n$. 
Observe that
\begin{equation}  \label{mu0}
\mu_{snN}=\sum\nolimits_{j=1}^N \I_{\{\eta_j=s\}}
\end{equation}
is the number of boxes containing $s$ balls.
Here $\I_{\{. \}}$ denotes  the indicator of a set.

In \cite{ChuFaz2010} an exponential inequality  for the tail of the conditional expectation of sums of centered independent random variables was obtained.
That inequality was applied  to prove analogues of the Law of the Iterated Logarithm and the  Strong Law of Large Numbers for conditional expectations. 
As corollaries, certain strong theorems for the generalized allocation scheme were obtained.

In the present paper we introduce the multi-colour version of the generalized allocation scheme.
Then we prove analogues of the Law of the Iterated Logarithm and Law of Large Numbers for the numbers of boxes containing fixed numbers of balls.
Our theorems are extensions of the results of \cite{ChuFaz2010}.
To prove our theorems, we use some basic inequalities of \cite{ChuFaz2010}.
\section{The multi-colour scheme and the main results}
\label{sectRes} \setcounter{equation}{0} 
%
Now we introduce the multi-colour version of the generalized allocation scheme. 
Let $K$ be a fixed positive integer, it will denote the number of different colours.
$n_i$ denotes the number of balls of $i$th colour, $i=1,2,\dots, K$.
Then $n= n_1 + \dots + n_K$ is the total number of balls.
Let $N$ denote the number of boxes.
The generalized allocation scheme of allocation $n$ balls into $N$ boxes when there are $K$ different colours among the balls is the following.

Let $\xi_{i1},\dots, \xi_{iN}$, $1\le i\le K$, be an array of independent non-negative integer valued random variables. 
Denote by $\etab_{1}, \dots, \etab_{N}$  a set of $K$-dimensional random vectors with the following property.
The  joint  distribution of $\etab_{1}, \dots, \etab_{N}$ admits representation
\begin{equation}  \label{EtaXi}
\P\{{\etab_{1}=\kb_{1}, \dots, \etab_{N}=\kb_{N}}\}  = \prod_{i=1}^K \P\{\xi_{i1}=k_{i1}, \dots, \xi_{iN}=k_{iN}\,\, |\,\, \xi_{i1}+\dots +\xi_{iN}=n_i\},
\end{equation}
for arbitrary non-negative integers valued vectors  ${\kb_j}=(k_{1j}, k_{2j},\dots k_{Kj})$, $1\le j\le N$ with $k_{i1}+\dots + k_{iN} = n_i$, $i=1, \dots ,K$.
Then we say that $\eta_{i1}, \dots, \eta_{iN}$, $i=1, \dots ,K$, obey a multi-colour generalized allocation scheme of placing of $n_i$ $i$th colour balls ($1\le i\le K$) into $N$ boxes. 
$\eta_{ij} = k_{ij}$ means that after allocating $n$ balls into $N$ boxes, in the $j$th box there are $k_{ij}$ balls of $i$th colour.


The number of boxes containing $s_i$ balls of $i$th colour for any  $1\le i\le K$, is
\begin{equation}  \label{mu}
\mu_{{\sb \nb}N}=  \sum_{j=1}^N \I_{\{{\etab_{j}=\sb}\}} =\sum_{j=1}^N\prod_{i=1}^K \I_{\{\eta_{ij}=s_i\}}
\end{equation}
where ${\sb}=(s_1, s_2,\dots, s_K)$ and   ${\nb}=(n_1, n_2,\dots ,n_K)$ .

Let $i$ be fixed.
Let the random variables $\xi_{i1},\dots, \xi_{iN}$ be identically distributed. 
In Kolchin's generalized allocation scheme usually power law distributions are considered (see [2]).
For our model it means that the random variables $\xi_{i1}, \dots, \xi_{iN}$ are distributed as follows: 
$q_{ik}=\P\{\xi_{i1}=k\}= ({b_{ik}\theta^k})/({k!B_i(\theta)}) $
where $b_{i0}, b_{i1}, \dots $ is a certain sequence of non-negative numbers, and $ B_i(\theta)=\sum_{k=0}^{\infty} {b_{ik}\theta^k}/{k!} $ for each fixed $i\in \{ 1,2,\dots, K\}$.
\begin{exmp}   \label{Poiss1}
The usual (that is not generalized) multi-colour allocation. 
Let $\xi_{ij}$ have Poisson distribution, i.e. for each $i$ and $j$ let $\P (\xi_{ij}= k) = \frac{\lambda_i^k}{ k!} e^{-\lambda_i}$, $k=0,1,\dots$.
Then
$$
\prod_{i=1}^K \P\{\xi_{i1}=k_{i1}, \dots, \xi_{iN}=k_{iN}\,\,\, |\,\,\, \xi_{i1}+\dots +\xi_{iN}=n_i\} = 
\prod_{i=1}^K \frac{n_i!}{k_{i1}!\dots k_{iN}!} \left( \frac{1}{N} \right)^{n_i}
$$
if $k_{i1}+\dots +k_{iN}=n_i$.
That is for each fixed $i$ the vector $\{\eta_{i1}, \dots, \eta_{iN}\}$ has polynomial distribution.
Now $\{\eta_{i1}=k_{i1}, \dots, \eta_{iN}=k_{iN}\} $ means that the box contents are $k_{i1}, \dots, k_{iN}$ after allocating $n_i$ balls into $N$ boxes 
during the usual random allocation procedure.
Moreover, the allocations of balls of different colours are independent.
We see that $\mu_{{\sb \nb}N}$ is the number of boxes containing $s_i$ balls of $i$th colour ($i=1,2,\dots, K$).

We see that the parameter $\theta$ used in the generalized allocation scheme is the same as $\lambda$ in the above usual allocation procedure.

We remark that in \cite{ChuFaz2005} using direct methods, we proved the following.
If $n,N \to \infty$ so that $n/N \to \lambda_0$, then
\begin{equation}  \label{limMu}
\frac{\mu_{rnN}}{N} \to  \frac{\lambda_0^r}{ r!} e^{-\lambda_0}
\end{equation}
almost surely.
Here $\mu_{rnN}$ is the number of boxes with $r$ balls (after allocating $n$ balls into $N$ boxes in the usual random allocation procedure).
In \cite{ChuFaz2009} we extended this result for generalized allocations.
 \hfill $\Box$
\end{exmp}
In what follows,  we consider  sequences of non-negative numbers $b_{i0}, b_{i1}, \dots$ with $b_{i0}>0, b_{i1}>0$ 
and assume that the convergence radius $R_i$ of the series
\begin{equation}   \label{B}
B_i(\theta)=\sum\nolimits_{k=0}^{\infty}\frac{b_{ik}\theta^k}{k!}
\end{equation}
is positive for each fixed $i$.
 Let us introduce the integer-valued random variable $\xi_i(\theta)$  (where $\theta >0$) with distribution
\begin{equation}   \label{xi}
 {\P}\{\xi_i(\theta)=k\}=\frac{b_{ik}\theta^k}{k!B_i(\theta)},\,\,\,k=0,1,2,\dots.
\end{equation}
 By \cite{Kol99}, one has
$$
m_i(\theta)=  {\E}\xi_i(\theta)=  \frac{\theta B'_i(\theta)}{B_i(\theta)}
$$
and
\begin{equation}  \label{sigma}
\sigma^2_i(\theta)=  {\D}^2\xi_i(\theta)=
\frac{\theta^2 B^{''}_i (\theta)}{B_i(\theta)}+ \frac{\theta B'_i(\theta)}{B_i(\theta)}- \frac{\theta^2 (B'_i(\theta))^2}{(B_i(\theta))^2}.
\end{equation}
The last equality implies that
\begin{equation}  \label{sigma_theta}
\sigma^2_i(\theta)=\theta m'_i(\theta).
\end{equation}
Let $0<\theta'_i<\theta^{''}_i<R_i$. 
If $\sigma^2_i(\theta)=0$ for some $\theta\in [\theta'_i, \theta^{''}_i]$, then the random variable $\xi_i(\theta)$ is a constant. 
Since $b_{i0}>0, b_{i1}>0$, the random variable $\xi_i(\theta)$ is not a constant.
Therefore $\sigma^2_i(\theta)$, $\theta\in [\theta'_i, \theta^{''}_i]$ is a positive continuous function. 
Consequently,
\begin{equation}  \label{C1C2}
0<C_1 =\inf_{\theta\in [\theta'_i, \theta^{''}_i]} \sigma_i^2(\theta) \le \sup_{\theta\in [\theta'_i, \theta^{''}_i]}\sigma_i^2(\theta)= C_2<\infty, \,\,\,\,1\le i\le K.
\end{equation}
By \eqref{sigma_theta} and \eqref{C1C2}, we have
$$
0<\frac{C_1}{\theta^{''}_i}=\inf_{\theta\in [\theta'_i, \theta^{''}_i]}m'_i(\theta) \le \sup_{\theta\in [\theta'_i, \theta^{''}_i]}m'_i(\theta)= \frac{C_2}{\theta'_i}<\infty,
\,\,\,\,1\le i\le K.
$$
So $m_i(\theta)$, $\theta\in [\theta'_i, \theta^{''}_i]$, is a positive, continuous, strictly increasing function. 
We will denote by $m^{-1}_i$ the inverse function of $m_i$.

We see that the random variable $\xi_i(\theta)$ has all moments, if $\theta<R_i$.

Throughout the paper let  $\xi_{i1}(\theta), \dots , \xi_{iN}(\theta)$ be independent copies of $\xi_i(\theta)$ where
$\xi_i(\theta)$ has distribution \eqref{xi}. 

We shall prove limit theorems when the number of boxes and the number of balls converge to infinity (the number of colours is fixed).

First we study the asymptotic behaviour along a fixed subsequence.
The subsequence will be indexed by $t$ ($t$ can be considered as the discrete time $t=1,2,\dots$, and $t$ will go to $\infty$).
Let $n_{it}$ denote the number of balls of $i$th colour at time $t$.
Let ${\nb}_t=(n_{1t},\dots n_{Kt})$.
Let $N_t$ be the number of boxes at time $t$,  assume that  $N_t<N_{t+1}$, $t\in\N$. 
Let ${\alphab}_t=(\alpha_{1t},\dots, \alpha_{Kt})=(n_{1t}/{N_t},\dots, n_{Kt}/N_t)$ be the ratios of the numbers of balls and the number of boxes, 
and let $\theta_{it}= m_i^{-1}(\alpha_{it})$,
\begin{equation}  \label{must}
p_{{\sb}t}=\prod_{i=1}^K\frac{b_{is_i}\theta^{s_i}_{it}}{s_i!B_i(\theta_{it})} ,\quad
\sigma^2_{{\sb}t}=p_{{\sb}t}(1-p_{{\sb}t})\quad \mbox{\text and} \qquad
\mu_{{\sb}t}=\sum_{i=1}^{N_t} \I_{\{\etab_i={\sb}\}}.
\end{equation}
At time $t$, $\mu_{{\sb}t}$ is the number of boxes containing $s_1,s_2,\dots s_K$ balls of colours $1,2,\dots, K$, respectively.

We remark that if $b_{is_i}=0$ for some $1\le i\le K$, then $\mu_{{\sb}t}=0$. 
So we can concentrate on the case $b_{is_i}>0$ for all $1\le i\le K$.

We will prove the following analogue of the Law of the Iterated Logarithm for $\mu_{{\sb}t}$.
\begin{thm}    \label{Theorem1.1} 
Let $0<\alpha'_i<\alpha''_i$ be such that $m_i^{-1}(\alpha''_i)<R_i$, $1\le i\le K$. 
Suppose that $0<\alpha'_i<\alpha_{it}<\alpha''_i$  for all $t\in \N$,  $1\le i\le K$.
Let ${\sb}$ be fixed and we assume  that $b_{is_i}>0$ for each fixed $1\le i\le K$. 
Then we have
$$
\limsup_{t\to\infty}\frac{|\mu_{{\sb}t}-\E\mu_{{\sb}t}|}{\sqrt{N_t\ln(N_t)}\sigma_{{\sb}t} } \le
4\sqrt{1+\frac{K}{2}}\,\,\,\, \,\,\,\mbox{\text almost surely.}
$$
\end{thm}
Now, we turn to a version of the previous theorem.
We shall consider the asymptotic behaviour in a sector instead of along a subsequence.
More precisely, the number of boxes and the numbers of balls will tend to infinity, while their ratios remain bounded.
Let $\alpha_i={n_i}/{N}$ be the ratio of the number of balls of $i$th colour and the number of boxes, 
and let $\theta_{i\alpha_i}= m_i^{-1}(\alpha_{i})$, $i=1,2,\dots, K$.
We shall use the notation
\begin{equation}  \label{psnN}
p_{{\sb}\nb N}=\prod_{i=1}^K\frac{b_{is_i}\theta^{s_i}_{i\alpha_i}}{s_i!B_i(\theta_{i\alpha_i})} ,\quad
\sigma^2_{{\sb}\nb N}=p_{{\sb}\nb N}(1-p_{{\sb}\nb N}).
\end{equation}

We can again concentrate on the case $b_{is_i}>0$ for all $1\le i\le K$.
\begin{thm}    \label{Theorem1.2} 
Let $0<\alpha'_i<\alpha''_i$ be such that $m_i^{-1}(\alpha''_i)<R_i$, $1\le i\le K$. 
Let ${\sb}$ be fixed and assume that $b_{is_i}>0$, $1\le i\le K$. 
Then for $\mu_{{\sb \nb}N}$ defined in \eqref{mu}, we have
\begin{equation}  \label{eqno{(1.10)}}
\limsup_{\begin{subarray}{c}{\nb\to\infty}, N\to\infty,\\ \alpha'_i<\alpha_i<\alpha''_i,\,\,1\le i\le K\end{subarray}}
\frac{|\mu_{{\sb \nb}N}-\E\mu_{{\sb}\nb N}|}{\sqrt{N\ln(N)}\sigma_{{\sb \nb}N}}  \le
4\sqrt{1+\frac{3K}{2}}\,\,\,\, \,\,\,\mbox{\text almost surely.}
\end{equation}
\end{thm}
\begin{cor}  \label{Corollary1.1}
Let $0<\alpha'_i<\alpha''_i$ be such that $m_i^{-1}(\alpha''_i)<R_i$, $1\le i\le K$. 
Let ${\sb}$ be fixed. 
Assume that $b_{i s_i} > 0$, $1\le i\le K$. 
Let $d_N\in\R$, $N\in \N $ be such that $\lim_{N\to\infty}\frac{\sqrt{N\ln(N)}}{d_N}=0$. 
Then for $\mu_{{\sb \nb}N}$ defined in \eqref{mu} we have
$$
\lim_{\begin{subarray}{c}{\bf n},N\to\infty,\\
\alpha'_i<\alpha_i<\alpha''_i,\,\,1\le i\le K\end{subarray}}
\frac{\mu_{{\sb \nb}N}-\E\mu_{{\sb \nb}N}}{d_N }=0 \,\,\,\, \,\,\,\mbox{\text almost surely.}
$$
\end{cor}
Using Corollary  \ref{Corollary1.1}, we obtain the following SLLN.
\begin{thm}    \label{Theorem1.3} 
Let $N\to\infty$,  $n_i\to \infty$ such that $n_i/N = \alpha_i\to\alpha_{i0}$, \,  $m_i(\theta_{\alpha_{i0}})=\alpha_{i0}$ and
$0<\theta_{\alpha_{i0}}<R_i$, $1\le i\le K$. 
Then for $\mu_{{\sb}\nb N}$ defined in \eqref{mu}, we have
\begin{equation}  \label{eqno{(1.11)}}
\frac{1}{N}\mu_{{\sb \nb}N}  \to  \prod_{i=1}^K\frac{b_{is_i}\theta^{s_i}_{\alpha_{i0}}}{s_i!B_i(\theta_{\alpha_{i0}})}
 \end{equation} 
almost surely as $N\to\infty$, $n_i\to\infty$ such that $n_i/N = \alpha_i\to\alpha_{i0}$, \, $1\le i\le K$.
\end{thm}
\begin{exmp}   \label{Poiss2}
We apply Theorem \ref{Theorem1.3} to the usual (that is not generalized) multi-colour allocation mentioned in Example \ref{Poiss1}. 
Let $\xi_{ij}$ have Poisson distribution with parameter $\lambda$.
Then $b_{ij}=1$, $B_i(\lambda)=e^{\lambda}$. 
As $m_i (\lambda)= \E \xi_{i1} = \lambda$, the parameter $n_i/N =\alpha_i$ coincides with $\lambda_i$.
Therefore, by \eqref{eqno{(1.11)}},
$$
\frac{1}{N}\mu_{{\sb \nb}N}  \to  \prod_{i=1}^K \frac{\lambda^{s_i}_{{i0}}}{s_i!e^{\lambda_{i0}}},
$$
as $N\to\infty$,  $n_i\to \infty$ such that $n_i/N = \lambda_{i0}$.
 \hfill $\Box$
\end{exmp}
\begin{exmp}   \label{GeomPoiss}
Consider the allocation of $n_1+n_2$ balls into $N$ boxes with the following method.
First we allocate $n_1$ indistinguishable balls into $N$ boxes.
Then we allocate $n_1$ distinguishable balls into the same $N$ boxes.
Assume that the two allocations are independent.
So we have the following allocation scheme
 $$
 \P\{{\etab_{1}=\kb_{1}, \dots, \etab_{N}=\kb_{N}}\}  = \prod_{i=1}^2 \P\{\xi_{i1}=k_{i1}, \dots, \xi_{iN}=k_{iN}\,\, |\,\, \xi_{i1}+\dots +\xi_{iN}=n_i\}.
$$
We know, that the first allocation is the same as the partition of the positive integer number $n$ into $N$ positive integer summands.
Therefore it is easy to prove (see \cite{Kol99}, Example 1.2.2) that in the generalized allocation scheme
$$
\P\{\eta_{\rm 1 1}=k_{11}, \dots, \eta_{\rm 1 N}=k_{1N}\}  =
\P\{\xi_{11}=k_{11}, \dots, \xi_{1N}=k_{1N}\,\,\, |\,\,\, \xi_{11}+\dots +\xi_{1N}=n_1\},
$$
we should choose geometrical random variables $ \xi_{1j}$, $j=1,\dots, N$.
That is with the usual notation $\P\{\xi_{1j}=k\}=p^k(1-p)$, $k=0,1,2\dots$, $0<p<1$.
Therefore $b_{1k}=k!$, $B_1(p) = 1/(1-p)$, $m_1(p) = \E \xi_{1j} = 1/p$.
So the relation of our parameters is $\theta = 1/\alpha$.
The second allocation scheme
$$
\P\{\eta_{\rm 2 1}=k_{1}, \dots, \eta_{\rm 2 N}=k_{N}\}  =
\P\{\xi_{21}=k_{21}, \dots, \xi_{2N}=k_{iN}\,\,\,|\,\,\, \xi_{21}+\dots +\xi_{2N}=n_2\},
$$
is the usual one with Poisson random variables $ \xi_{2j}$.
That is $\P\{\xi_{2j}=k\}=e^{-\lambda}\frac{\lambda^k}{k!}$, $k=0,1,2\dots$, $0<\lambda$.
Therefore, by Theorem \ref{Theorem1.3}, we have
$$
\frac{1}{N}\mu_{{\sb\nb}N} \to \left(\frac{1}{\alpha_{10}}\right)^{s_1}\left(1-\frac{1}{\alpha_{10}}\right) e^{-\alpha_{20}}\frac{\alpha_{20}^{s_2}}{s_2!}
,\,\,\,\mbox{\text as}\,\,\, n_i,N\to\infty\,\,\,\mbox{\text such that}\,\,\,\,n_i/N= \alpha_i\to\alpha_{i0}, \,\,i=1,2,
$$
almost surely, where $1<  \alpha_{10}$, $0< \alpha_{20}$.
\hfill $\Box$
 \end{exmp}
\begin{rem}    \label{Remark1}
For each fixed $i\in\{1,2,\dots,K\}$ let $\zeta_{ij}$, $j\in{\N}$, be independent identically distributed random variables with finite expectation $m_i= \E\zeta_{ij}$.
Assume that  $\{ \zeta_{ij}\, : \, j\in{\N}\}$, $1\le i\le K$, are independent sets of random variables. 
By the Kolmogorov Law of Large Numbers we have
 $$
 \frac{1}{N}\sum_{j=1}^N\prod_{i=1}^K\zeta_{ij}\to\prod_{i=1}^K m_{i},\,\,\,\,\,\mbox{\rm almost surely} \,\,\,\,\,\mbox{\rm as}\,\,\,\,\,\,N\to\infty.
 $$
Theorem \ref{Theorem1.3} is an analogue of this assertion for an independent sequence of sets of dependent indicators.
 \hfill $\Box$
\end{rem} 
Now return to scheme \eqref{EtaXi}.
Let $s$ be a fixed non-negative integer number.
Denote by $\mu_{s{\nb}N}$ the number of boxes containing altogether $s$ balls regardless of their colour.
Then we have
\begin{equation}  \label{eqno{(1.12)}}
\mu_{s{\nb}N}=\sum_{{\sb}=(s_1, s_2,\dots s_K), s_1 + s_2 + \dots +s_K=s} \mu_{{\sb \nb}N}
=\sum_{s_1 + s_2 + \dots +s_K=s}\sum_{j=1}^N\prod_{i=1}^K\I_{\{\eta_{ij}=s_i\}}.
\end{equation}
Since the number of vectors ${\sb}=(s_1, s_2,\dots ,s_K)$ such that $s_1 + s_2 + \dots +s_K=s$ is finite, Theorem \ref{Theorem1.3} implies the following.
\begin{cor}  \label{Corollary1.2}
Let $N\to\infty$,  $n_i\to \infty$ such that $n_i/N =\alpha_i\to\alpha_{i0}$ and $m_i(\theta_{\alpha_{i0}})=\alpha_{i0}$ for some $0<\theta_{\alpha_{i0}}<R_i$, $1\le i\le K$. 
Then for $\mu_{s{\nb}N}$ defined in \eqref{eqno{(1.12)}}, we have
\begin{equation}  \label{eqno{(1.13)}}
\frac{1}{N}\mu_{s{\nb}N} \to \sum_{ s_1 + s_2 + \dots +s_K=s} \prod_{i=1}^K\frac{b_{is_i}\theta^{s_i}_{\alpha_{i0}}}{s_i!B_i(\theta_{\alpha_{i0}})}
\end{equation}
almost surely as $N\to\infty$,  $ n_i\to\infty$ such that $n_i/N=\alpha_i\to\alpha_{i0}$, $1\le i\le K$.
\end{cor}
\begin{rem}    \label{Remark2}
Consider again the usual allocation.
That is the case then $\xi_{ij}$ has  Poisson distribution with parameter $\lambda$. 
Then $b_{ij}=1$, $B_i(\lambda)=e^{\lambda}$, $m_i(\lambda)=\lambda$. 
Therefore the limit in \eqref{eqno{(1.13)}} has the following form
$$
\sum_{s_1 + s_2 + \dots +s_K=s} \prod_{i=1}^K \frac{\lambda^{s_i}_{i0}}{s_i!e^{\lambda_{i0}}}=
e^{-(\lambda_{10}+\lambda_{20}+\dots +\lambda_{K0})}\frac{(\lambda_{10}+\lambda_{20}+\dots +\lambda_{K0})^s}{s!}.
$$
This limit coincides with the limit of $\frac{1}{N}\mu_{snN}$ at the usual random allocation of $n$ balls into $N$ boxes
when $\frac{n}{N}\to \lambda_{10}+\lambda_{20}+\dots +\lambda_{K0}$  (see \eqref{limMu}). 

So the following problem arises: what kind of distributions of $\xi_{ij}=\xi_{ij}(\theta)$ satisfy this property. 
More precisely: which power series distributions satisfy that for each $K$ the multi-colour limit in \eqref{eqno{(1.13)}}
is the same as the single colour limit when $\frac{n}{N}\to \alpha_{10}+\alpha_{20}+\dots +\alpha_{K0}$.

If it is valid, then for $B_i=B$, $s=0$ and $\alpha_{10}=\alpha_{20}=\dots =\alpha_{K0}=\frac{\alpha}{K}$ using \eqref{eqno{(1.13)}}, we obtain
\begin{equation}  \label{eqno{(1.15)}}
 \left(\frac{b_0}{B\left(\frac{\alpha}{K}\right)}\right)^K  =  \frac{b_0}{B(\alpha)}.
\end{equation}
 Since
  $$
 \left(\frac{b_0}{B\left(\frac{\alpha}{K}\right)}\right)^K  \to  e^{-(b_1/b_0)\alpha}, \quad \mbox{\rm as}\,\,\,\,K\to\infty,
 $$
 from \eqref{eqno{(1.15)}} it follows that $B(\alpha)=b_0 e^{(b_1/b_0)\alpha}$, therefore $\xi_{ij}$ has a  Poisson distribution.
  \hfill $\Box$
 \end{rem}
\begin{rem}    \label{Remark3}
Now we show that in scheme \eqref{EtaXi0} the distribution of $\eta_1,\dots , \eta_N$ determines the underlying power series distribution (up to a multiplicative constant).
So consider scheme \eqref{EtaXi0} both with $\P(\xi_i=k) = \frac{b_k\theta^k}{k!B(\theta)}$ and $P(\xi_i=k) =\frac{\tilde{b}_k\theta^k}{k!\tilde{B}(\theta)}$. 
Let $n,N\to\infty$ such that $\frac{n}{N}\to\alpha$, where $\alpha<\infty$.
 By Theorem \ref{Theorem1.3},
$\frac{1}{N}\mu_{0nN}\to\frac{b_0}{B(\alpha)}$ almost surely and at the same time $\frac{1}{N}\mu_{0nN}\to\frac{\tilde{b}_0}{\tilde{B}(\alpha)}$ almost surely 
as $n,N\to\infty$ such that $n/N\to \alpha$.
So $\frac{b_0}{B(\alpha)} = \frac{\tilde{b}_0}{\tilde{B}(\alpha)}$ which implies $B(\alpha) = (b_0/\tilde{b}_0)\tilde{B}(\alpha)$.
 \hfill $\Box$
 \end{rem}
\section{Proofs}
\label{sectProofs} \setcounter{equation}{0} 
In order to prove Theorem \ref{Theorem1.1} we need the following lemma. 
Consider the events
\begin{equation}  \label{eqno{(2.1)}}
A_{iN}(n_i)=\{\xi_{i1}(\theta_{i\alpha_i})+\xi_{i2}(\theta_{i\alpha_i})+\dots + \xi_{iN}(\theta_{i\alpha_i})=n_i\}
\end{equation}
where $\theta_{i\alpha_i}=m^{-1}_i(\alpha_i)$, $\alpha_i = n_i/N$. 
\begin{lem} \label{Lemma2.1}
(Lemma 3.1 in \cite{ChuFaz2009}.)
Let $0< \alpha'_i< \alpha''_i$ be such that $m^{-1}_i( \alpha''_i)<R_i$.  
Then there exists $N_0\in{\N}$ with the following property: 
if $n_i, N$ are positive integers such that $N>N_0$  and $ \alpha'_i\le\alpha_i\le \alpha''_i$,
then we have
\begin{equation}  \label{eqno{(2.2)}}
 \P (A_{iN}(n_i)) > \frac{1}{4\sigma_i(\theta_{i\alpha_i})\sqrt{N}}.
\end{equation}
\end{lem}
{\sc Proof of Theorem \ref{Theorem1.1}.} 
We want to apply Theorem 3.1 of \cite{ChuFaz2010}.
By the independence of the random variables $\xi_{ij}$, $1\le i\le K$, $1\le j\le N$, and by \eqref{EtaXi}, we see that the general framework of that theorem covers our model.


Let $A_t=\cap_{i=1}^K A_{iN}(n_{it})$.
Using first \eqref{EtaXi} and \eqref{eqno{(2.2)}}, then \eqref{C1C2}, we obtain
$$
\P(A_t)> \frac{1}{(4\sqrt{N_t})^K}\prod_{i=1}^K\frac{1}{\sigma_i(\theta_{it})} \ge\frac{C''}{N_t^{K/2}}.
$$
As $C''>0$, condition (3.1) of \cite{ChuFaz2010} is valid with $\beta=\frac{K}{2}$. 
We have $\sigma^2_{{\sb}t}= p_{{\sb}t}(1-p_{{\sb}t})$. 
Here
$$
p_{{\sb}t}=\prod_{i=1}^K \frac{b_{is_i}\theta^{s_i}_{it}}{s_i!B_i(\theta_{it})}>
\prod_{i=1}^K\frac{b_{is_i}(\theta'_i)^{s_i}}{s_i!B_i(\theta''_{i})}>0
$$
as $b_{is_i}>0$. 
Since $b_{i0}$, $b_{i1}>0$, therefore we have
$$
1-p_{{\sb}t}>\prod_{i=1}^K\left(1-\frac{b_{is_i}\theta^{s_i}_{it}}{s_i!B_i(\theta_{it})}\right)>c>0.
$$
 So for some $c_2>0$ the relation $\sigma^2_{{\sb}t}\ge c_2$ is valid for each $t$, i.e. (3.2) of \cite{ChuFaz2010} is satisfied. 
 Thus Theorem \ref{Theorem1.1} follows from Theorem 3.1 of \cite{ChuFaz2010}.
 \hfill $\Box$

 In order to prove Theorem \ref{Theorem1.2}, we need the following corollary of Theorem 2.1 of \cite{ChuFaz2010}. 
 We shall use the following notation.
 $\sigma^2$  is the variance of $\I_{\{(\xi_{1j},\dots\xi_{Kj})={\sb}\}}$, 
 $\rho$ is the variance of $\I_{\{(\xi_{1j},\dots\xi_{Kj})={\sb}\}}-\I'_{\{(\xi_{1j},\dots\xi_{Kj})={\sb}\}}$, 
 where $\I'_{\{(\xi_{1j},\dots\xi_{Kj})={\sb}\}}$ is an independent copy of $\I_{\{(\xi_{1j},\dots\xi_{Kj})={\sb}\}}$.
 $A_{{\nb}N}=\cap_{i=1}^K A_{iN}(n_{i})$ where $A_{iN}(n_{i})$ is defined in \eqref{eqno{(2.1)}}.
 $Z$ is a centered Gaussian random variable with variance $1$.
\begin{lem} \label{Lemma2.2}
Let $\varepsilon\ge 4\sqrt{2}\sigma$. 
Then we have
\begin{equation}  \label{eqno{(2.3)}}
\P\left\{ \frac{|\mu_{{\sb \nb}N}-\E \mu_{{\sb\nb}N}|}{\sqrt{N}}
\ge \varepsilon\right\} \le \frac{\sqrt{2}}{{\P}(A_{{\nb}N})} e^{-\frac{\varepsilon^2}{16\sigma^2}}(1+B),
\end{equation}
where
\begin{equation}  \label{eqno{(2.4)}}
B=B(N, \sigma,\varepsilon)=
\frac{\rho}{32}\left(\frac{\varepsilon^2}{8\sigma^4\sqrt{N}}\right)^2f_2\left(\frac{2\varepsilon^2}{8\sigma^4\sqrt{N}}\right)
+\O\left(\frac{8\sigma^2}{\varepsilon^2}\right)
\end{equation}
and
$$
f_2(x)=2\E\{Z^2\exp(x|Z|)\}-1.
$$
\end{lem}
{\sc Proof of Theorem \ref{Theorem1.2}.}
First we remark that limsup exists and it is unique for a multiindex sequence in a sector.
That is 
$$
\limsup_{\begin{subarray}{c}{\nb\to\infty}, N\to\infty,\\ \alpha'_i<\alpha_i<\alpha''_i,\,\,1\le i\le K\end{subarray}}
\frac{|\mu_{{\sb \nb}N}-\E\mu_{{\sb}\nb N}|}{\sqrt{N\ln(N)}\sigma_{{\sb \nb}N}}
$$
is well-defined.

Since  $b_{is_i}>0$, $b_{i0}>0$, $b_{i1}>0$ and $\alpha'_i < \alpha < \alpha''_i $,   for all $1\le i\le K$, 
therefore we have $0<c_1 \le \sigma^2_{{\sb \nb}N}\le 1/4$, and $\rho \sigma_{{\sb \nb}N} \rho_{{\sb \nb}N}\le 1/2$. 
Let $z$ be a fixed positive number.
Therefore for  $B$ from \eqref{eqno{(2.4)}}, we have $B(N, \sigma, z \ln(N) \sigma_{{\sb \nb}N})  \le L<\infty$ 
as $N,{\nb}\to \infty$ such that $\alpha'_i<\alpha_i<\alpha''_i$,  $1\le i\le K$.

Now let $z>4\sqrt{1+\frac{3K}{2}}$. 
Then $\frac{z^2}{16}-\frac{K}{2}-K>1$. 
Therefore, by \eqref{eqno{(2.3)}} and \eqref{eqno{(2.2)}}, we have
$$
\sum_{N=N_0+1}^{\infty}\sum_{\begin{subarray}{c}
N\alpha'_i<n_i<N\alpha''_i,\\1\le i\le K\end{subarray}}
\P\left\{ \frac{|\mu_{{\sb \nb}N}-\E \mu_{{\sb \nb}N}|}{\sqrt{N\ln(N)}\sigma_{{\sb \nb}N}} \ge z\right\}=
$$
$$
= \sum_{N=N_0+1}^{\infty}\sum_{\begin{subarray}{c}
N\alpha'_i<n_i<N\alpha''_i,\\1\le i\le K\end{subarray}}
\P\left\{ \frac{|\mu_{{\sb\nb}N}-\E \mu_{{\sb\nb}N}|}{\sqrt{N}} \ge z\sqrt{\ln(N)}\sigma_{{\sb\nb}N}\right\}\le
$$
$$
\le \sup_{\begin{subarray}{c}N>N_0,\\\,\,
N\alpha'_i<n_i<N\alpha''_i,\\1\le i\le K\end{subarray}}
(1+B)\sum_{N=N_0+1}^{\infty}\sqrt{2}\frac{N^{K/2}}{C''}\prod_{i=1}^K(N(\alpha''_i-\alpha'_i))e^{-\frac{\ln(N) z^2}{16}}\le
$$
$$
\le \sup_{\begin{subarray}{c}N>N_0,\\\,\,
N\alpha'_i<n_i<N\alpha''_i,\\1\le i\le K\end{subarray}}
(1+B)\frac{\sqrt{2}}{C''} \left(\prod_{i=1}^K((\alpha''_i-\alpha'_i))\right)\sum_{N=N_0+1}^{\infty}N^{-\frac{z^2}{16}+K+\frac{K}{2}}<\infty.
$$
Consequently, by the Borel-Cantelli lemma, we have
$$
\limsup_{\begin{subarray}{c}{\nb},N\to\infty,\\
\alpha'_i<\alpha_i<\alpha''_i,\,\,1\le i\le K\end{subarray}}
\frac{|\mu_{{\sb \nb}N}-\E\mu_{{\sb}\nb N}|}{\sqrt{N\ln(N)}\sigma_{{\sb \nb}N}}  \le z\,\,\,\, \,\,\,\mbox{\text almost surely.}
$$
This and the choice of $z$ imply \eqref{eqno{(1.10)}}. 
 \hfill $\Box$

In order to prove Theorem \ref{Theorem1.3}, we need the following lemma.
\begin{lem}  \label{Lemma2.3}
(Lemma 4.1 of \cite{ChuFaz2009}.)
Let $i\in\{ 1, \dots , N \}$ be fixed.
Let $0< \alpha'_i< \alpha^{''}_i$ be such that $m^{-1}_i( \alpha^{''}_i)<R_i$.
Let $\theta_{\alpha_i} = m^{-1}(\alpha_i)$.
 For any $0\le s_i<\infty$, as $N, n_i\to \infty$, uniformly for $\alpha'_i<\alpha_i<\alpha''_i$ we have
$$
{\E}\left(\frac{1}{N}\mu_{s_in_iN}\right) =\P\{\eta_{ij}=s_i\}=
\frac{b_{is_i}\theta_{\alpha_i}^{s_i}}{s_i!B_i(\theta_{\alpha_i})}
\cdot\sqrt{\frac{N}{N-1}}\cdot 
\frac{ \frac{1}{\sqrt{2\pi}}\exp\left\{-\frac{\left(n_i-s_i-\frac{n_i}{N}(N-1)\right)^2}{2
\sigma^2_i(\theta_{\alpha_i})(N-1)}\right\}+\o_{1i}(1)}
{\frac{1}{\sqrt{2\pi}}+\o_{2i}(1) },
$$
where $\o_{1i}, \o_{2i}\to 0$.
Here $\mu_{s_in_iN}$ is the number of boxes containing $s_i$ balls after placing $n_i$ balls into $N$ boxes.
\end{lem}
{\sc Proof of Theorem \ref{Theorem1.3}.} 
Observe that
\begin{equation}  \label{eqno{(2.5)}}
\frac{1}{N} \mu_{{\sb\nb}N}= \frac{\mu_{{\sb\nb}N}-\E\mu_{{\sb\nb}N}}{N}+\E\left(\frac{1}{N}\mu_{{\sb\nb}N}\right).
\end{equation}
By  \eqref{mu} and Lemma  \ref{Lemma2.3}, we obtain
$$
{\E}\left(\frac{1}{N}\mu_{{\sb\nb}N}\right)
=\prod_{i=1}^K\P\{\eta_{ij}=s_i\}=\prod_{i=1}^K \frac{b_{is_i}\theta_{\alpha_i}^{s_i}}{s_i!B_i(\theta_{\alpha_i})}
\cdot\sqrt{\frac{N}{N-1}}
\cdot \frac{\frac{1}{\sqrt{2\pi}}\exp\left\{-\frac{\left(n_i-s_i-\frac{n_i}{N}(N-1)\right)^2}{2\sigma^2_i(\theta_{\alpha_i})(N-1)}\right\}+\o_{1i}(1)}
{\frac{1}{\sqrt{2\pi}}+\o_{2i}(1) }.
$$
Therefore we have
\begin{equation}  \label{eqno{(2.6)}}
\E\left(\frac{1}{N}\mu_{{\sb\nb}N}\right)\to
\prod_{i=1}^K\frac{b_{is_i}\theta^{s_i}_{\alpha_{i0}}}{s_i!B_i(\theta_{\alpha_{i0}})},\,\,\,\mbox{\text as}\,\,\, n_i,N\to\infty\,\,\,
\mbox{\text such that}\,\,\,\,n_i/N= \alpha_i\to\alpha_{i0}, \,\,\,\,1\le i\le K. 
\end{equation}
By Corollary 1.1
\begin{equation}  \label{eqno{(2.7)}}
\frac{\mu_{{\sb\nb}N}-\E \mu_{{\sb\nb}N}}{N}\to 0,\,\,\,\mbox{\text as}\,\,\, n_i,N\to\infty\,\,\,\mbox{\text such that}
\,\,\,\,n_i/N= \alpha_i\to\alpha_{i0}, \,\,\,\,1\le i\le K,\,\,\,\mbox{\text almost surely}. 
\end{equation}
From  \eqref{eqno{(2.5)}}, by \eqref{eqno{(2.6)}} and \eqref{eqno{(2.7)}}, follows \eqref{eqno{(1.11)}}. 
 \hfill $\Box$
\end{document}